\documentclass{segabs}

\usepackage{latexsym,amssymb,amsmath,amsthm, amsfonts, float}
\usepackage{enumerate}

\def\RR{{\mathbb R}}

\def\cal{\mathcal}

\def\cO{{\cal O}}

\usepackage{array}
\newcolumntype{L}{>{\displaystyle} l}

\usepackage{color}

\renewcommand{\v}{\mathbf{v}}
\renewcommand{\H}{\mathbf{H}}
\renewcommand{\u}{\mathbf{u}}
\newcommand{\f}{\mathbf{f}}

\newcommand{\x}{\mathbf{x}}

\newcommand{\beq}{\begin{equation}}
\newcommand{\eeq}{\end{equation}}

\newcommand{\bit}{\begin{itemize}}
\newcommand{\eit}{\end{itemize}}
\newcommand{\ben}{\begin{enumerate}}
\newcommand{\een}{\end{enumerate}}

\title{A short note on the nested-sweep polarized traces method for the 2D Helmholtz equation}

\author{Leonardo Zepeda-N\'u\~nez\textsuperscript{(*)}, and Laurent Demanet\\
Department of Mathematics and Earth Resources Laboratory, Massachusetts Institute of Technology }

\lefthead{Zepeda-N\'u\~nez \& Demanet}
\righthead{A short note on the nested-sweep polarized traces method for the 2D Helmholtz equation}

\begin{document}

\maketitle

\begin{abstract}

We present a variant of the solver in \cite{ZepedaDemanet:the_method_of_polarized_traces}, for the 2D high-frequency Helmholtz equation in heterogeneous acoustic media. By changing the domain decomposition from a layered to a grid-like partition, this variant yields improved asymptotic online and offline runtimes and a lower memory footprint. The solver has online parallel complexity that scales \emph{sublinearly} as $\cO \left( \frac{N}{P} \right)$, where $N$ is the number of volume unknowns, and $P$ is the number of processors, provided that $P = \cO(N^{1/5})$. The variant in \cite{ZepedaDemanet:the_method_of_polarized_traces} only afforded $P = \cO(N^{1/8})$. Algorithmic scalability is a prime requirement for wave simulation in regimes of interest for geophysical imaging.

\end{abstract}

\section{Introduction}

Define the Helmholtz equation, in a bounded domain $\Omega \subset \RR^2$, with frequency $\omega$ and squared slowness $m(\x) = 1/c^2(\x)$, by
\begin{equation} \label{eq:Helmholtz}
	\left( - \triangle - m(\x) \omega^2 \right) u(\x) = f(\x)
\end{equation}
with absorbing boundary conditions. In this note, Eq \ref{eq:Helmholtz} is discretized with a 5-point stencil, and the absorbing boundary conditions are implemented via a perfectly matched layer (PML) as described by \cite{Berenger:PML}. This leads to a linear system of the form $\H \u = \f$. Let $N$ be the total number of unknowns of the linear system and $n = N^{1/2}$ the number of points per dimension. In addition, suppose that the number of point within the PML grows as $\log(N)$. There is an important distinction between:
\begin{itemize}
\item the offline stage, which consists of any precomputation involving $\H$, but not $\f$; and
\item the online stage, which involves solving $\H\u = \f$ for many different right-hand sides $\f$.
\end{itemize}

By online complexity, we mean the runtime for solving the system once in the online stage. The distinction is important in situations like geophysical wave propagation, where offline precomputations are often amortized over the large number of system solves with the same matrix $\H$.

\cite{ZepedaDemanet:the_method_of_polarized_traces} proposed a hybrid between a direct and an iterative method, in which the domain is decomposed in $L$ layers. Using the Green's representation formula to couple the subdomains, the problem is reduced to a boundary integral system at the layer interfaces, solved iteratively, and which can be preconditioned efficiently using a polarizing condition in the form of incomplete Green's integrals. Finally, the integral system and the preconditioner are compressed, yielding a fast application.

\begin{figure}[H]
\centering
\includegraphics[trim= 40mm 90mm 40mm 10mm, angle =90  ,clip,  width=7.8cm]{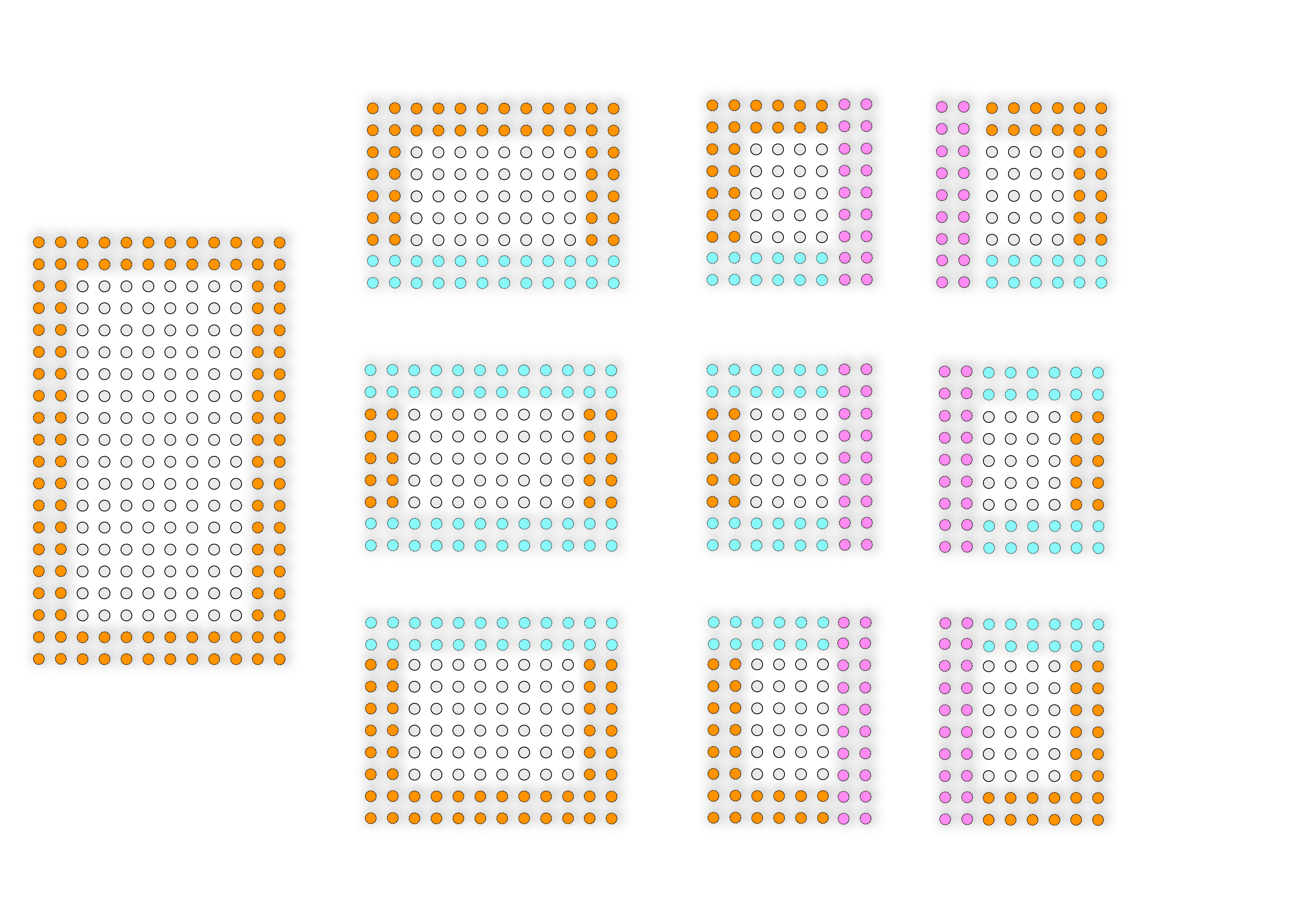}
\caption{ Nested Decomposition in cells. The orange grid-points represent the PML for the original problem, the light-blue represent the artificial PML between layers, and the pink grid-points represent the artificial PML between cells in the same layer.} \label{fig:DDM}
\vspace{-.3cm}
\end{figure}

Let $P$ for the number of nodes in a distributed memory environment. In \cite{ZepedaDemanet:the_method_of_polarized_traces}, each layer is associated with one node ($P = L$), and the method's online runtime is $\cO(N/P)$ as long as $P = \cO(N^{1/8})$. The method presented in this note instead uses a nested domain decomposition, with a first decomposition into $L \sim \sqrt{P}$ layers, and a further decomposition of each layer into $L_c \sim \sqrt{P}$ cells, as shown in Fig. \ref{fig:DDM}. The resulting online runtime is now $\cO(N/P)$ as long as $P = \cO(N^{1/5})$.

The main drawback of the method of polarized traces in \cite{ZepedaDemanet:the_method_of_polarized_traces} is its offline precomputation that involves computing and storing interface-to-interface local Green's functions. In 3D this approach would become impractical given the sheer size of the resulting matrices. To alleviate this issue, the nested-sweep polarized traces method presented in this paper involves an equivalent \emph{matrix-free} approach that relies on local solves with sources at the interfaces between layers. Given the iterative nature of the preconditioner, solving the local problems naively would incur a deterioration of the online complexity. This deterioration can be circumvented if we solve the local problems inside the layer via the same boundary integral strategy as earlier, in a \emph{nested} fashion. This procedure can be written as a factorization of Green's integrals in block-sparse factors.

Finally, the offline complexity is much reduced; instead of computing large Green's functions for each layer, we compute much smaller interface-to-interface operators between the interfaces of adjacent cells within each layer, resulting in a lower memory requirement.

% In particular, we argue that for fixed number $P$ of nodes, a grid-like decomposition in $P$ cells ($\cO(P^{1/2})$ layers, with $\cO(P^{1/2})$ cells inside each layer) yields lower storage and lower runtimes than a decomposition in $L=P$ layers. We can show that for $P \lesssim N^{1/5}$ the method becomes strictly sublinear, whereas the method of polarized traces is strictly sublinear only when $P \lesssim N^{1/8}$.

% If the frequency scales as $\omega \sim \sqrt{n}$, the regime in which second order finite-differences are expected to be accurate, we obtain, empirically, a storage requirement of $\cO(P^{3/8}N^{5/8}+ P\log(N)^{5/4} )$, instead of $\cO(N^{5/8}P)$ for the layered case. The cost of one iteration is given by $\cO(P^{3/8} N^{5/8}+ P\log(N)^{5/4} )$ in the grid-like case, in contrast with $\cO(P N^{5/8})$ in the layered case. This allows us to obtain a better scaling of $P$ with respect to the number of unknowns. Indeed, the grid-like decomposition yields an empirical online complexity of $\cO(N^{8/11})$ with $P \sim N^{3/11}$, compared to $\cO(N^{13/16})$ with $P \sim N^{3/16}$ for the layered case.

\subsection{Literature Review}

%Solving the time-harmonic wave equation with a scalable algorithm in the high frequency regime is a question of great interest for applications in geophysical imaging, and has attracted great attention during the last years.

 Classical iterative methods require a large number of iterations to converge and standard algebraic preconditioners often fail to improve the convergence rate to tolerable levels (see \cite{Erlangga:Helmholtz}). Classical domain decomposition methods also fail because of internal reverberations (see \cite{Gander:why_is_difficult_to_solve_helmholtz_problems_with_classical_iterative_methods}). Only recently, the development of a new type of precondtioners have shown that domain decomposition with accurate transmission boundary conditions is the right mix of ideas for simulating propagating high-frequency waves in a heterogeneous medium. To a great extent, the approach can be traced back to the AILU method of \cite{GanderNataf:ailu_for_hemholtz_problems_a_new_preconditioner_based_on_an_analytic_factorization}. The first linear complexity claim was perhaps made in the work of \cite{EngquistYing:Sweeping_H}, followed by \cite{CStolk_rapidily_converging_domain_decomposition}.
Other authors have since then proposed related methods with similar properties, including \cite{Chen_Xiang:a_source_transfer_ddm_for_helmholtz_equations_in_unbounded_domain} and \cite{GeuzaineVion:double_sweep}; however, they are often difficult to parallelize, and they usually rely on distributed linear algebra such as \cite{Poulson_Engquist:a_parallel_sweeping_preconditioner_for_heteregeneous_3d_helmholtz} and \cite{Tsuji_Poulson:sweeping_preconditioners_for_elastic_wave_propagation}, or highly tuned multigrid methods such as \cite{Bhowmik_Stolk:A_multigrid_method_for_the_helmholtz_equation_with_optimized_coarse_grid_corrections} and \cite{Calandra_Grattonn:an_improved_two_grid_preconditioner_for_the_solution_of_3d_Helmholtz}. As we put the final touches to this note, \cite{Liu_Ying:Recursive_sweeping_preconditioner_for_the_3d_helmholtz_equation} proposed a recursive sweeping algorithm closely related to the one presented in this note.

In a different direction, much progress has been made on making direct methods efficient for the Helmholtz equation. Such is the case the HSS solver in \cite{Wang:H_multifrontal} and hierarchical solver in \cite{Gillman_Barnett_Martinsson:A_spectrally_accurate_solution_technique_for_frequency_domain_scattering_problems_with_variable_media}. The main issue with direct methods is the lack of scalability to very large-scale problems due to the memory requirements.

To our knowledge, \cite{ZepedaDemanet:the_method_of_polarized_traces} and this note are the only two references to date that make claims of sublinear $O(N/P)$ online runtime for the 2D Helmholtz equation. We are unaware that this claim has been made in 3D yet, but we are confident that the same ideas have an important role to play there as well.

\section{Methods}
\subsection{Polarized Traces}
%We call an algorithm is scalable if it has lower asymptotic runtime in the parallel case than in the best sequential case.

This section reviews the method in \cite{ZepedaDemanet:the_method_of_polarized_traces}. Provided that the domain is decomposed in $L$ layers, it is possible to reduce the discrete Helmholtz equation to a discrete integral system at the interfaces between layers, using the Green's representation formula. This leads to an equivalent discrete integral system of the form
\begin{align} \label{eq:discrete_integral_system}
\underline{\mathbf{M}} \, \underline{\u} =\underline{\f}.
\end{align}
The online stage has three parts: first, at each layer, a local solve is performed to build $\underline{\f}$ in Eq. \ref{eq:discrete_integral_system}; then Eq. \ref{eq:discrete_integral_system} is solved obtaining the traces $\underline{\u}$ of the solution at the interfaces betwen layers; finally, using the traces, a local solve is performed at each layer to obtain the solution $\u$ in the volume. The only non-parallel step is solving Eq. \ref{eq:discrete_integral_system}. Fortunately, there is a physically meaningful way of efficiently preconditioning this system.

\begin{figure}[H]
\centering
\includegraphics[trim = 20mm 22mm 16mm 17mm, clip, width=3cm]{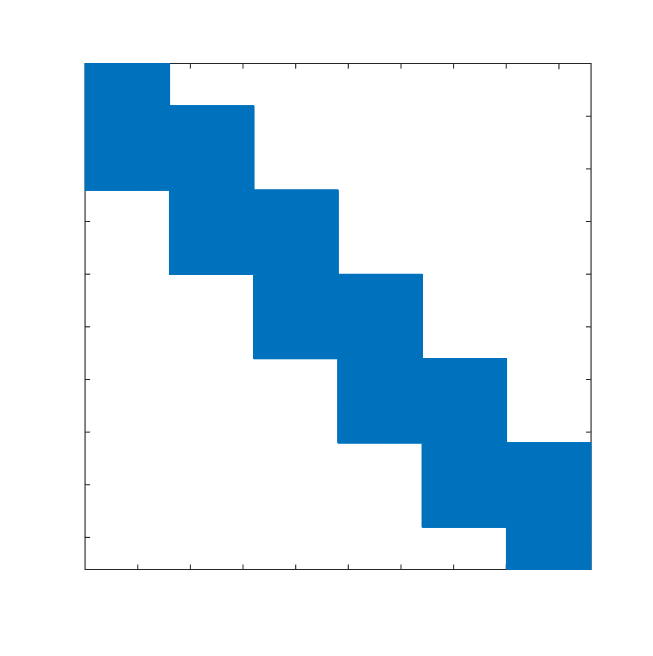}
\includegraphics[trim = 20mm 22mm 16mm 17mm, clip, width=3cm]{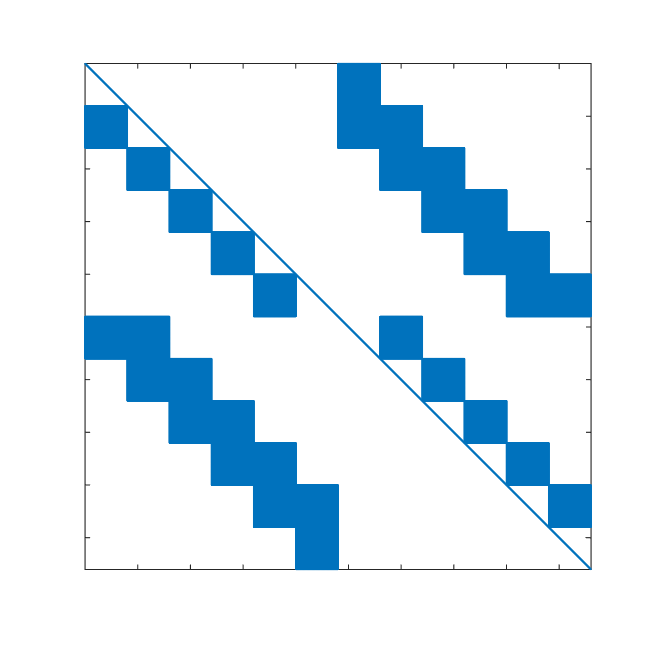}
\caption{Sparsity pattern of  $\underline{\mathbf{M}}$ (left) and $\underline{\underline{\mathbf{M}}}$ (right).} \label{fig:M_polarized}
\vspace{-.3cm}
\end{figure}

The key is to write yet another equivalent problem where local Green's functions are used to perform polarization into one-way components. A wave is said to be polarized as up-going at an interface $\Gamma$ when
\begin{align*}
0  =  \int_{\Gamma}\partial_{\nu_{\x'}} G(\x, \x')  u^{\uparrow}(\x') dx'- \int_{\Gamma} G(\x, \x') \partial_{\nu_{\x'}} u^{\uparrow}(\x') dx'   ,
\end{align*}
as long as $x$ is below $\Gamma$, and vice-versa for down-going waves. These polarization conditions create cancellations within the discrete integral system, resulting in an easily preconditionable system for the polarized interface traces $\underline{\u}^{\uparrow}$ and $\underline{\u}^{\downarrow}$ such that $\underline{\u} =\underline{\u}^{\uparrow} + \underline{\u}^{\downarrow}$.
The resulting polarized system is
\begin{align}
\underline{\underline{\mathbf{M}}} \underline{\underline{\u}} =\underline{\underline{\f}}, \qquad  \underline{\underline{\u}} = \left( \begin{array}{c}
\underline{\u}^{\downarrow} \\
\underline{\u}^{\downarrow}
\end{array} \right). \label{eq:integral_polarized}
\end{align}
The matrix $\underline{\underline{\mathbf{M}}}$ has the structure depicted in Fig. \ref{fig:M_polarized}. In particular, we write
\begin{equation}
\underline{\underline{\mathbf{M}}} = \left [ \begin{array}{cc}
\underline{\mathbf{D}}^{\downarrow} & \underline{\mathbf{U}} \\
\underline{\mathbf{L}}              &  \underline{\mathbf{D}}^{\uparrow}
\end{array}  \right ],
\end{equation}
% As in \cite{ZepedaDemanet:the_method_of_polarized_traces} we can easily define a block Jacobi preconditioner,
% \begin{equation}
% P^{\text{Jacobi}} \left( \begin{array}{c} \v^{\downarrow} \\ \v^{\uparrow}
% \end{array} \right ) =  \left( \begin{array}{c} (\underline{\mathbf{D}}^{\downarrow})^{-1} \v^{\downarrow} \\ (\underline{\mathbf{D}}^{\uparrow})^{-1} \v^{\uparrow}
% \end{array} \right ),
% \end{equation}
where $ \underline{\mathbf{D}}^{\downarrow}$ and $\underline{\mathbf{D}}^{\uparrow}$ are, respectively, block lower and upper triangular matrices with identity diagonal blocks, thus trivial to invert via block backsubstitution.
In this note, we use the Gauss-Seidel Preconditioner,
\begin{equation}
P^{\text{GS}} 	\left( 	\begin{array}{c} \v^{\downarrow} \\
									  \v^{\uparrow}
					 	\end{array}
				\right) =  	\left(
								\begin{array}{c}(\underline{\mathbf{D}}^{\downarrow})^{-1} \v^{\downarrow} \\
												(\underline{\mathbf{D}}^{\uparrow})^{-1} \left ( \v^{\uparrow} -\underline{\mathbf{L}}  (\underline{\mathbf{D}}^{\downarrow})^{-1} \v^{\downarrow} \right)
								\end{array}
							\right ).
\end{equation}
The system in Eq. \ref{eq:integral_polarized} is solved using GMRES preconditioned with $P^{\text{GS}}$.
Moreover, one can use an adaptive partitioned low-rank (PLR, also known as $\mathcal{H}$-matrix) fast algorithm for the application of integral kernels, in expressions such as the one above. The implementation details are in \cite{ZepedaDemanet:the_method_of_polarized_traces}.

\subsection{Nested Solvers}

% The main caveat of the method of polarized traces is its offline runtime. Although the offline complexity is lower than direct methods, the constants are much larger. The main bottleneck is the computation of the interface-to-interface Green's functions used to form the polarized system  $\underline{\underline{\mathbf{M}}}$. Even though the computation is embarrassingly parallel, using $\cO(nP)$ nodes is often impractical.  A more reasonable assumption would be to have $P$ nodes, which implies an offline runtime of $\cO(N^{3/2}/P)$ which can be prohibitively expensive. The other issue related to the computation of the Green's functions is its memory footprint, which can be intractable for the 3D version of this algorithm.

This section now describes the new solver.

\begin{figure}[H]
\vspace{-.3cm}
\centering
\includegraphics[trim= 20mm 14mm 10mm 30mm,clip,  width=7cm]{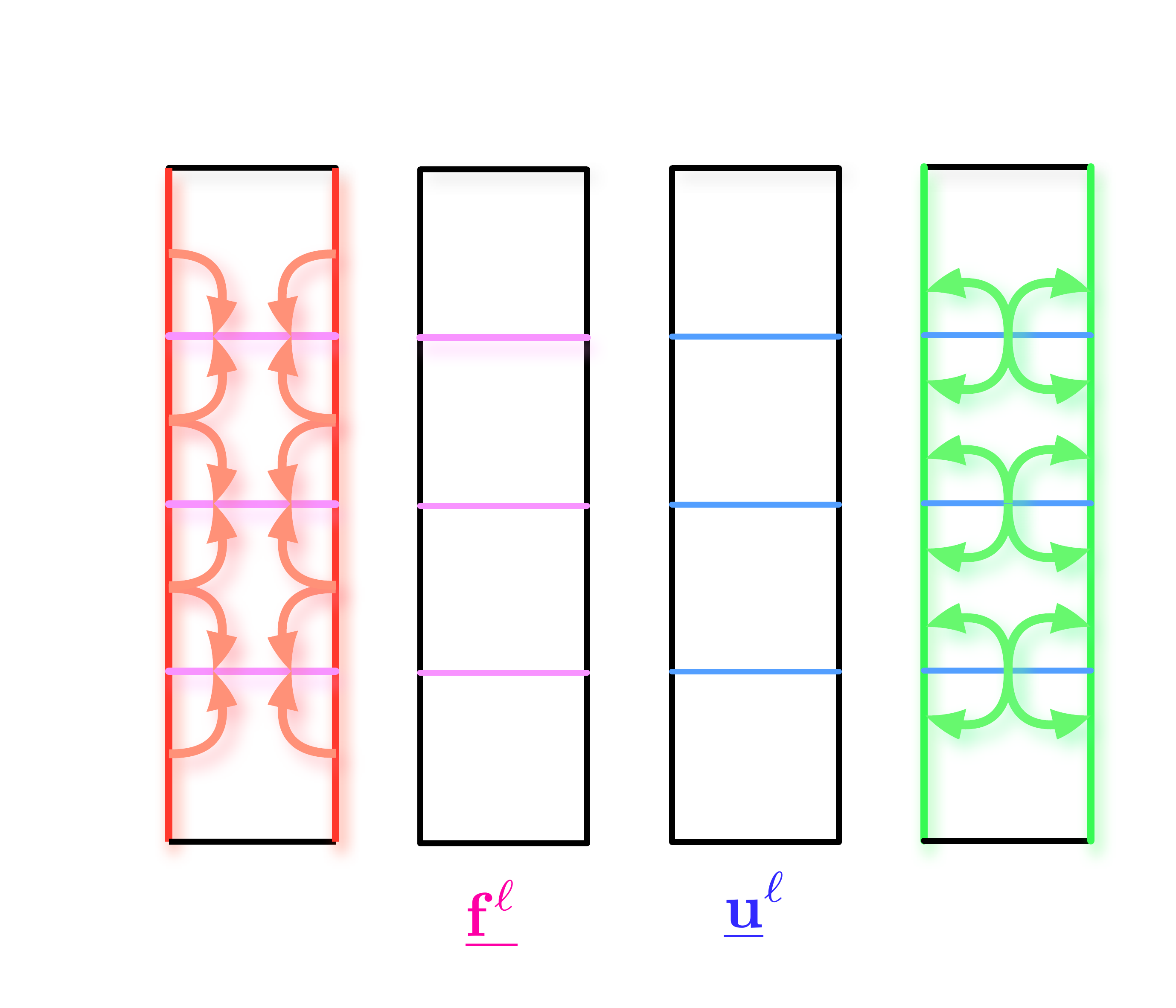}
\caption{ Sketch of the application of the Green's functions. The sources are in red (left) and the sampled field in green (right). The application uses the inner boundaries as proxies to perform the solve.} \label{fig:fact_sketch}
\vspace{-.3cm}
\end{figure}

One key observation is that the polarized matrix $\underline{\underline{\mathbf{M}}}$ can be applied to a vector in a matrix-free fashion, addressing the offline bottleneck and memory footprint. Each block of $\underline{\underline{\mathbf{M}}}$ is a Green's integral, and its application to a vector is equivalent to sampling a wavefield generated by sources at the boundaries. The application of each Green's integral to a vector $\v$, in matrix-free form, consists of three steps: from $\v$ we form the sources at the interfaces (red interfaces in Fig. \ref{fig:fact_sketch} left), we perform a local direct solve inside the layer, and we sample the solution at the interfaces (green interfaces in Fig. \ref{fig:fact_sketch} right).

 % Indeed, parsing $\v$, forming the sources at the interfaces of each layer, solving the local problem at each interface and then sampling the resulting wavefields at the interfaces, is equivalent to an application of $\underline{\underline{\mathbf{M}}}$ to $\v$.

% The potential disadvantage of matrix-free approach lies on the application of the preconditioner, which is translated to $\cO(L)$ direct solves applied sequentially. This sequential operation would greatly deteriorate the online complexity if applied naively.

From the analysis of the rank of off-diagonal blocks of the Green's functions, we know that the Green's integrals can be compressed and applied in $\cO(n^{3/2})$ time, but this requires precomputation and storage of the Green's functions. The matrix-free approach does not need expensive precomputations, but would naively perform a direct local solve in the volume, resulting in an application of each Green's integral in $\cO(n^2/L)$ complexity  (assuming that a good direct method is used at each layer). This note's contribution is to propose ``nested" alternatives for the inner solves that result in a lower $\cO(L_c (n/L)^{3/2})$ complexity (up to logs) for the application of each Green's integral.

%This becomes problematic when applying the preconditioner, which involves $L$ sequential applications. There seems to be a trade-off between offline and online cost; however, we explain below how to reorganize the computation in order to improve both online and offline runtimes.

%We follow the matrix-free approach but instead of a direct solve we use a nested solver, i.e., we use the same reduction used in the whole $\Omega$ to each layer. We reduce the local problem at each layer to solving a discrete integral system analog to Eq. \ref{eq:discrete_integral_system} with a layered decomposition in the transverse direction, given by

Analogously to Eq. \ref{eq:discrete_integral_system}, denote the boundary integral reduction of the local problem within each layer as
\begin{align} \label{eq:local_integral_system}
\underline{\mathbf{M}}^{\ell} \underline{\u}^{\ell} =\underline{\f^{\ell}},\qquad  \text{for }\, \ell = 1,..,L_c;
\end{align}
where the unknowns $\underline{\u}^{\ell}$ are on the interfaces between the $L_c$ cells.

The nested solver uses the inner boundaries as proxies to perform the local solve inside the layer efficiently. The application of the Green's integral can be decomposed in three steps. Using precomputed Green's functions at each cell we evaluate the wavefield generated from the sources to form $\underline{\f^{\ell}}$ (from red to pink in Fig. \ref{fig:fact_sketch} left). This operation can be represented by a sparse block matrix $\mathbf{M}_f^{\ell}$.  Then we solve Eq. \ref{eq:local_integral_system} to obtain  $\underline{\u}^{\ell}$ (from pink to blue in Fig. \ref{fig:fact_sketch} right). Finally we use the Green's representation formula to sample the wavefield at the interfaces (from blue to green in Fig. \ref{fig:fact_sketch}), this operation is represented by another sparse-block matrix $\mathbf{M}_u^{\ell}$.

The algorithm described above leads to the factorization
\begin{equation} \label{eq:factorization}
G^{\ell} = \mathbf{M}_f^{\ell} \left( \underline{\mathbf{M}}^{\ell} \right )^{-1} \mathbf{M}_u^{\ell},
 \end{equation}
in which the blocks of $\mathbf{M}_f^{\ell}$  and $\mathbf{M}_u^{\ell}$ are dense, but highly compressible in PLR form.

% In theory, we would need to perform a local solve inside each cell to build $\underline{\f^{\ell}}$; however, given that the sources are supported on the boundaries between layers, we can precompute the Green's functions linking vertical and horizontal interfaces within each cell. This precomputation provides matrices of size order $n/L \times n/L_c $, which are highly compressible.

% To solve Eq. \ref{eq:local_integral_system} we used two different methods with comparable online complexities, but different constants, and different offline runtimes. Both approaches lead to a fast method to solve the system.

% Once the $ \underline{\u}^{\ell} $ has been computed, it is possible to sample the wavefield at the interfaces using the Green's representation formula and the Green's functions linking the boundaries of the cell. Once again, this procedure can be efficiently compressed, yielding a fast application.

\subsection{Nested polarized traces}

To efficiently apply the Green's integrals using Eq. \ref{eq:factorization}, we need to solve Eq. \ref{eq:local_integral_system} efficiently. A first possibility would be to use the method of polarized traces in a recursive fashion. The resulting recursive polarized traces method has empirically the same scalings as those found in \cite{ZepedaDemanet:the_method_of_polarized_traces} when the blocks are compressed in partitioned low rank (PLR) form. However, the constants are much larger.

\subsection{A compressed-block LU solver}
A better alternative is to use a compressed-block LU solver to solve Eq. \ref{eq:local_integral_system}. Given the banded structure of $\underline{\mathbf{M}}^{\ell}$ (see Fig. \ref{fig:M_polarized}), we perform a block LU decomposition without pivoting (without fill-in), which leads to the factorization
\begin{equation}
G^{\ell} = \mathbf{M}_f^{\ell} \left( \underline{\mathbf{U}}^{\ell} \right )^{-1} \left( \underline{\mathbf{L}}^{\ell} \right )^{-1} \mathbf{M}_u^{\ell}.
 \end{equation}

Furthermore, the diagonal blocks of the LU factors are inverted explicitly so that the forward and backward substitutions are reduced to a series of matrix-vector products. Finally, the blocks are compressed in PLR form, yielding a fast solve. This compressed direct strategy is not practical at the outer level because of its prohibitive offline cost, but is very advantageous at the inner level.

The online complexity is comparable with the complexity of the nested polarized traces, but with much smaller constants. The LU factorization has a $ \cO \left( N^{3/2}/P + \log(N)^{3/2} \right)$ offline cost; which is comparable to the polarized traces method, but with much smaller constants.

%This offline cost can in principle be further decreased by using $\cH$-matrices algebra as in \cite{Hackbusch:Hierarchical_matrices}.

\section{Complexity}

Table \ref{table:complexity} summarizes the complexities and number of processors at each stage for both methods. For simplicity we do not count the logarithmic factors from the nested dissection; however, we consider the logarithmic factors coming from the extra degrees of freedom in the PML.
\begin{table}
    \begin{center}
        \begin{tabular}{|c|c|c|}
            \hline
            Step        & $N_{\text{nodes}}$  & Complexity per node \\
            \hline
            LU factorizations   & $O(P)$        & $O\left( (N/P + \log(N))^{3/2} \right)$ \\
            Green's functions   & $O(P)$       & $O\left( (N/P + \log(N))^{3/2} \right) $  \\
            \hline
            Local solves            & $O(P) $   & $O\left( N/P + \log(N)^2 \right )$ \\
            Sweeps          & 1                 & $O( P (N/P  + \log(N)^2 )^{\alpha}) $ \\
            Recombination       & $O(P) $       & $O\left( N/P + \log(N)^2  \right) $ \\
            \hline
        \end{tabular}
    \end{center}
    \caption{Complexity of the different steps of the preconditioner, in which $\alpha$ depends on the compression of the local matrices, thus on the scaling of the frequency with respect to the number of unknowns. Typically $\alpha = 3/4$.}\label{table:complexity}
    \vspace{-.3cm}
\end{table}

If the frequency scales as $\omega \sim \sqrt{n}$, the regime in which second order finite-differences are expected to be accurate, we observe $\alpha = 5/8$; however, we assume the more conservative value $\alpha = 3/4$. The latter is in better agreement with a theoretical analysis of the rank of the off-diagonal blocks of the Green's functions. In such scenario we have that the blocks of $\mathbf{M}_u^{\ell}$ and $\mathbf{M}_f^{\ell}$ can be compressed in PLR form, resulting in a fast application in $\cO( L_c(n/L + \log(n))^{3/2} )$ time, easily parallelizable among $L_c$ nodes. Solving Eq. \ref{eq:local_integral_system} can be solved using either the direct compressed or the nested polarized traces in $\cO( L_c(n/L + \log(n))^{3/2} )$ time. This yields the previously advertised runtime of  $\cO( L_c(n/L + \log(n))^{3/2} )$ for each application of the Green's integral.

\vspace{-.1cm}
\begin{figure}[H]
\centering
\includegraphics[trim = 12mm 12mm 25mm 10mm, clip, width=7.8cm]{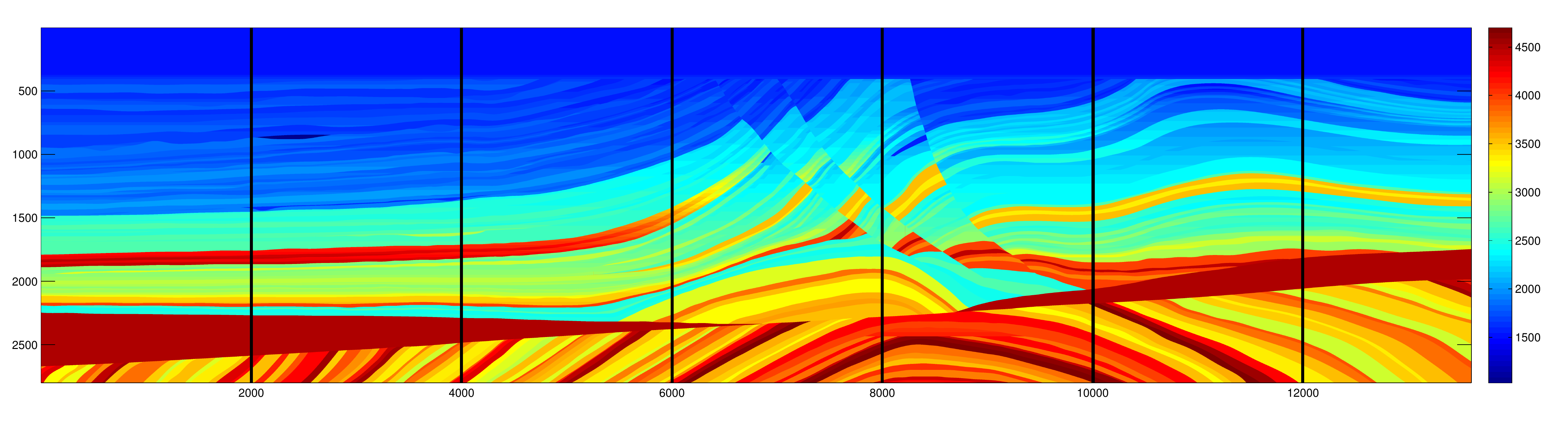}
\end{figure}

\vspace{-.6cm}

\begin{table}[H]
    \begin{center}
        \begin{tabular}{|c|c|r|r|r|r|}
            \hline
            $N$ 		& $f [Hz] $ & $10 \times 2$  		& $40 \times 8$ 	 & $100 \times 20$   \\
            \hline
            % $ 88 \times 425$  &      1.4 & \textbf{(3)}  0.85  & \textbf{(3)}   8.09  & \textbf{(3)} 67.6 \\
            % $175 \times 850$  &      2.7 & \textbf{(3)}  1.38  & \textbf{(3)}   10.2  & \textbf{(3)} 84.4  \\
            % $350 \times 1700$ &      5.5 & \textbf{(3)}  3.22  & \textbf{(3)}   16.7  & \textbf{(3)} 110  \\
            % $700 \times 3400$  &     11.2 & \textbf{(3)}  7.94 & \textbf{(3)}   30.1  & \textbf{(3)} 172 \\
$ 88 \times 425$  &      7.71  &	 \textbf{(3)} 0.89 &	\textbf{(3)} 15.6 & 	\textbf{(4)} 97.9 \\
$175 \times 850$  &      11.1  &	 \textbf{(3)} 1.48 &	\textbf{(3)} 17.7 & 	\textbf{(3)} 105 \\
$350 \times 1700$ &      15.9  &	 \textbf{(3)} 2.90 &	\textbf{(3)} 22.1 & 	\textbf{(4)} 106 \\
$700 \times 3400$  &     22.3  &	 \textbf{(3)} 5.58 &	\textbf{(3)} 31.3 & 	\textbf{(4)} 126 \\
$1400 \times 6800$	&    31.7  &  	\textbf{(3)} 10.5  &	\textbf{(3)} 47.9 & 	\textbf{(4)} 176 \\
            \hline
        \end{tabular}\caption{Number of GMRES iterations (bold) required to reduce the relative residual to $10^{-5}$, along with average execution time (in seconds) of one GMRES iteration using the compressed direct method, for different $N$ and $P = L \times L_c$. The frequency is scaled such that $f = \omega/\pi \sim \sqrt{n}$ and the sound speed is given by the Marmousi2 model (see \cite{Marmousi_2}).}   \label{table:numerical_results}
    \end{center}
    \vspace{-.3cm}
\end{table}

To apply the Gauss-Seidel preconditioner we need $\cO(L)$ applications of the Green's integral, resulting in a runtime of $ \cO( L\cdot L_c(n/L + \log(n))^{3/2} )$ to solve Eq. \ref{eq:integral_polarized}. Using the fact that $L \sim \sqrt{P}$ and $L_c \sim \sqrt{P}$ and adding the contribution of the other steps of the online stage; we have that the overall online runtime is given by $\cO(P^{1/4} N^{3/4}+ P\log(N)^{3/4} + N/P + \log(N)^2)$, which is sublinear and of the form $O(N/P)$ up to logs, provided that $P  = \cO(N^{1/5})$.

\begin{figure}[H]
\centering
\includegraphics[trim= 6.5cm 4.56cm 6.3cm 3.84cm, clip,  width=5.3cm]{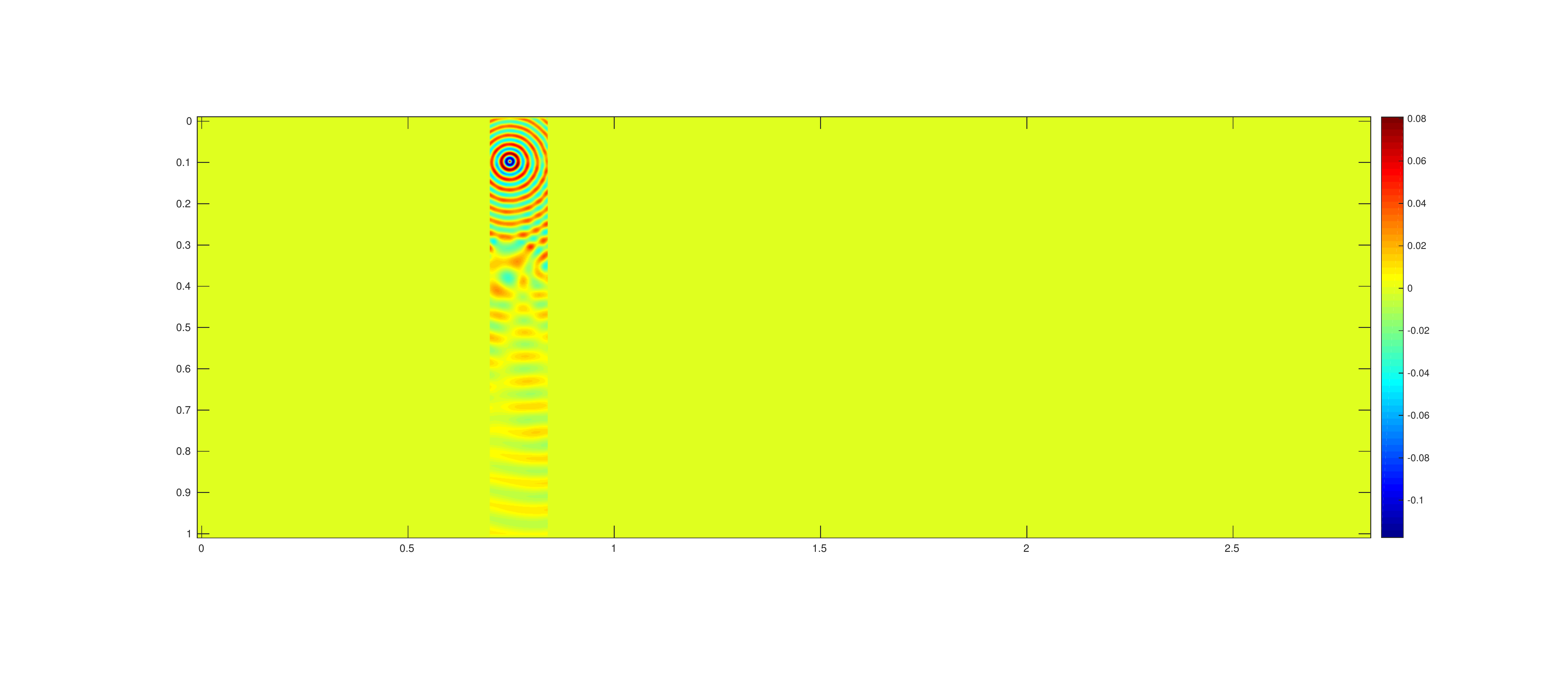}
\includegraphics[trim= 6.5cm 4.56cm 6.3cm 3.84cm, clip,  width=5.3cm]{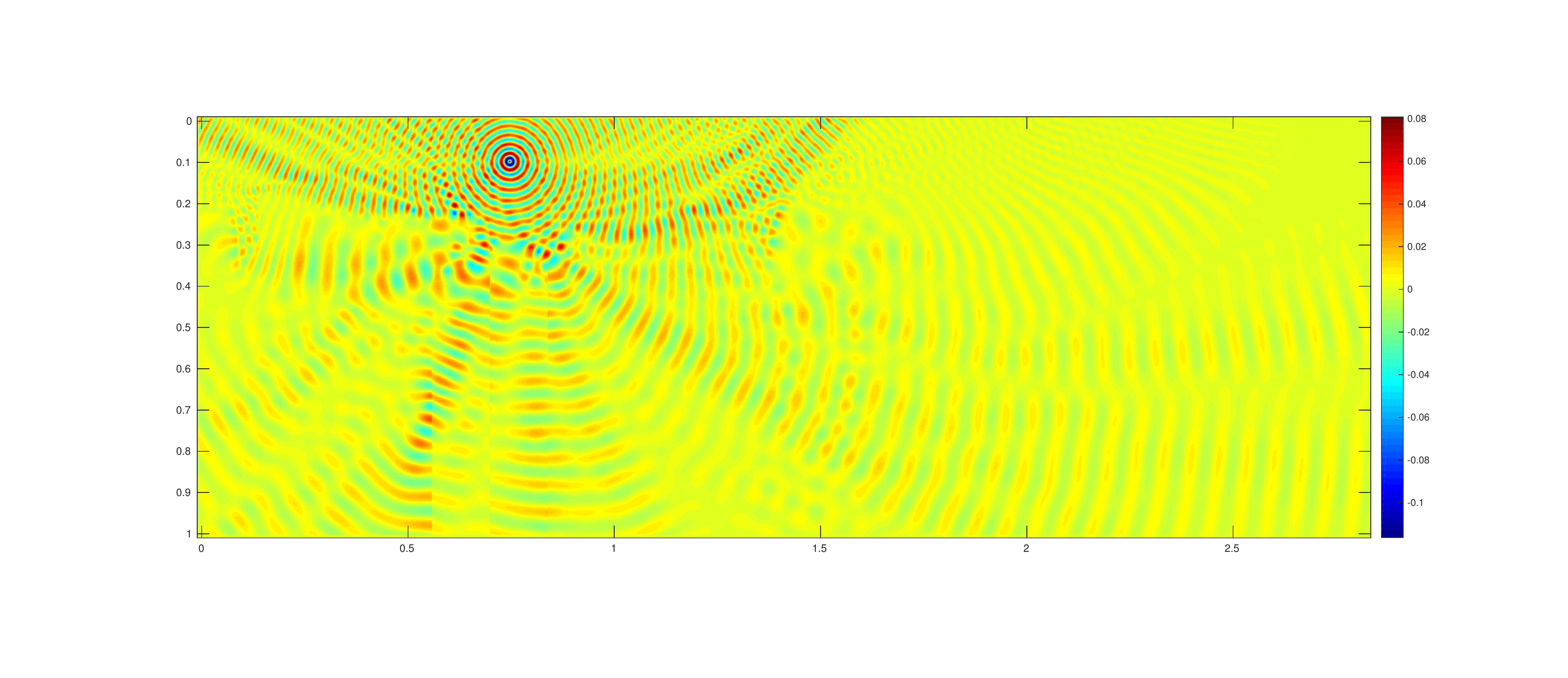}
\includegraphics[trim= 6.5cm 4.56cm 6.3cm 3.84cm, clip,  width=5.3cm]{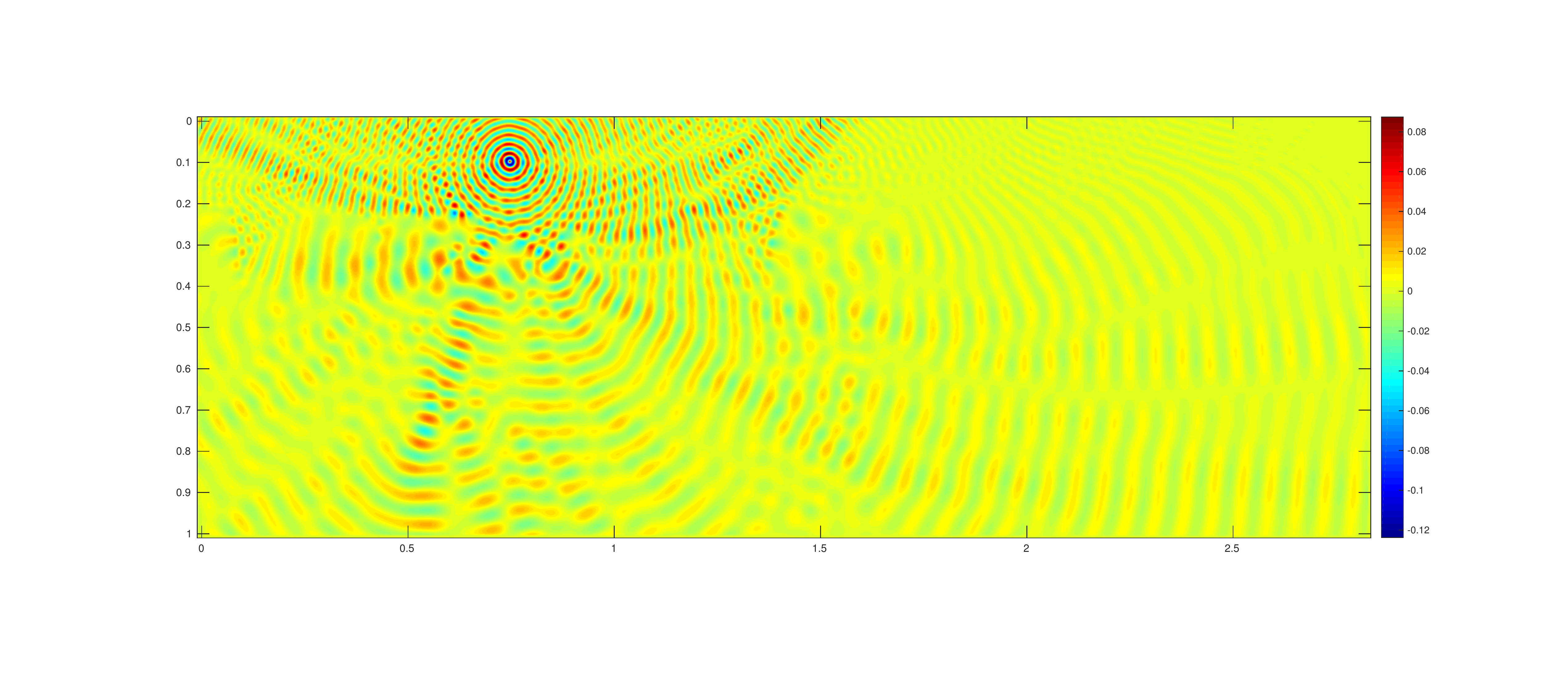}
\includegraphics[trim= 6.5cm 4.56cm 6.3cm 3.84cm, clip,  width=5.3cm]{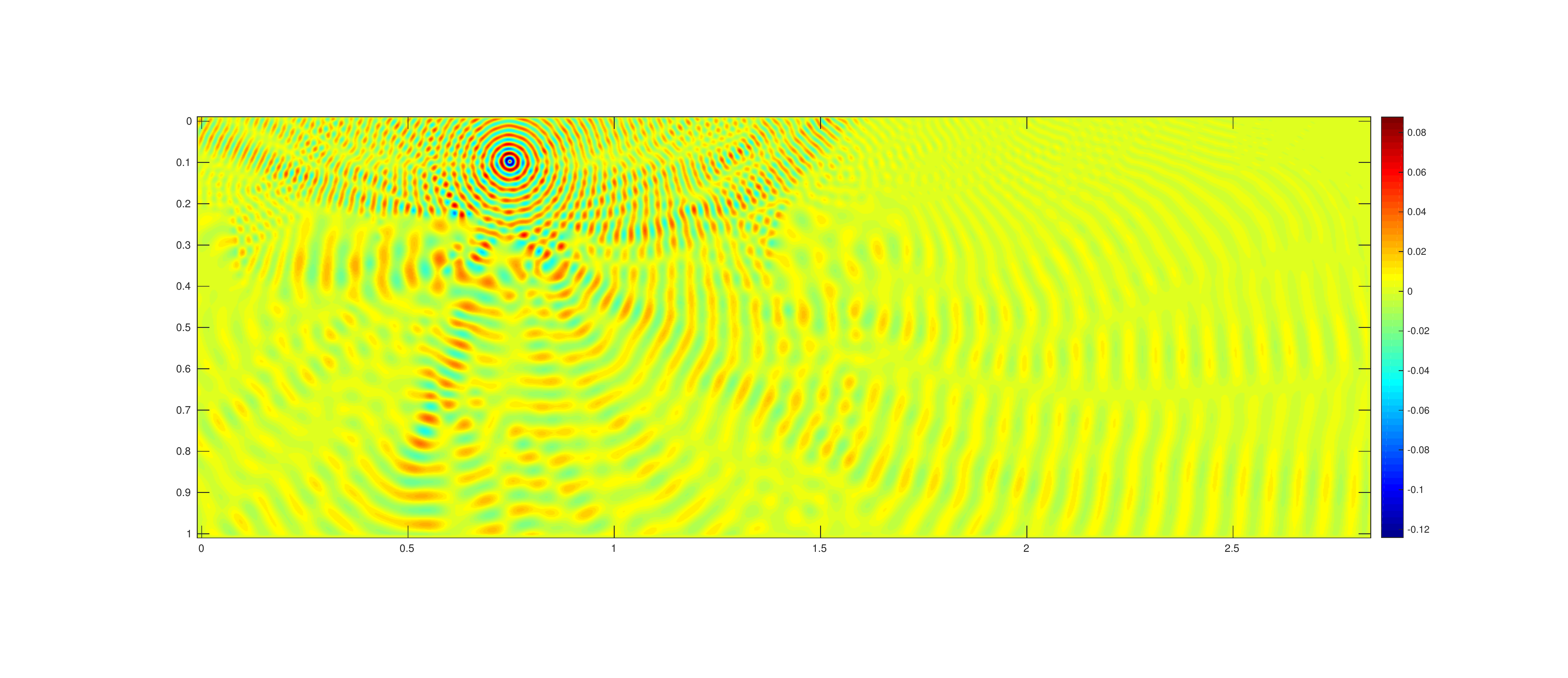}
\caption{ Two iteration of the preconditioner, from top to bottom: initial guess with only local solve; first iteration, second iteration, converged solution.} \label{fig:iteration}
\end{figure}
\vspace{-.3cm}
%  we obtain, empirically, a storage requirement of $\cO(P^{3/8}N^{5/8}+ P\log(N)^{5/4} )$ and a runtime of $\cO(P^{3/8} N^{5/8}+ P\log(N)^{5/4} )$ to solve the integral system. However, a theoretical analysis of the ranks of the off-diagonal block of the Green's functions yields a storage requirement of $\cO(P^{1/4}N^{3/4}+ P\log(N)^{3/4} )$  and a complexity  of $\cO(P^{1/4} N^{3/4}+ P\log(N)^{3/4})$ to solve the integral system (the mistmatch is due to discretization effects). These scalings lead to an overall runtime  of $\cO(P^{1/4} N^{3/4}+ P\log(N)^{3/4} + N/P + \log(N)^2)$, which is strictly sublinear if $L \lesssim N^{1/5}$. In the layered case we obtain a overall runtime of $\cO(P N^{3/4}+ N/P )$, which is strictily sublinear only then $P \lesssim N^{1/8}$.

% When the frequency $\omega \sim \sqrt{n}$, regime in which the second-order finite-differences discretization is expected to be accurate, we obtain an empirical scaling such that $\alpha = 5/8$ (theoretically we can only show $\alpha = 3/4$). When we considered the three stages of the online runtime, we obtain an overall runtime of $O( P (N/P  + \log(N)^2 )^{5/8} +  N/P + \log(N)^2 )$ , which is strictly sublinear if $P \lesssim N^{3/11}$ (theoretically we only have $P \lesssim N^{1/5}$).
Finally, the memory footprint is $\cO(P^{1/4} N^{3/4}+ P\log(N)^{3/4} + N/P + \log(N)^2)$ and the communication cost for the online part is $\cO(n\sqrt{P})$, which represents an asymptotic improvement with respect to \cite{ZepedaDemanet:the_method_of_polarized_traces}, in which the storage and communication cost are $\cO(PN^{3/4} + N/P + \log(N)^2)$ and $\cO(nP)$, respectively.

\section{Numerical results}

Fig. \ref{fig:iteration} depicts the fast convergence of the method. After a couple of iterations the exact and approximated solution are visually indistinguishable.

\begin{figure}[H]
\centering
\includegraphics[trim= 0mm 0mm 0mm 0mm,clip, width=7cm]{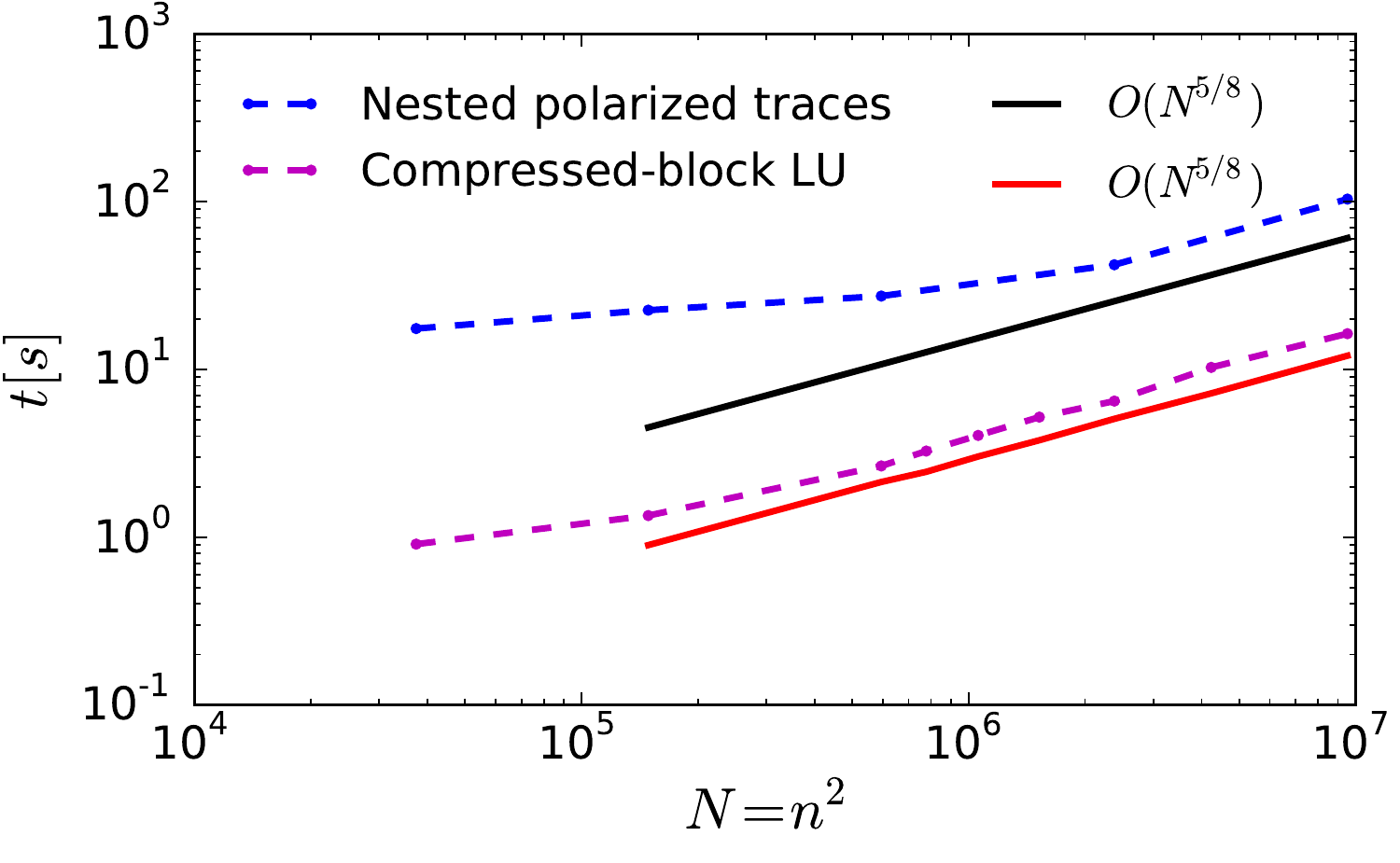}
\caption{ Runtime for one GMRES iteration using the two different nested solves, for $L=9$ and $L_c = 3$, and $\omega \sim \sqrt{n}$.} \label{graph:complexity}
\vspace{-.6cm}
\end{figure}

Table \ref{table:numerical_results} shows the sublinear scaling for one GMRES iteration, with respect to the degrees of freedom in the volume. We can observe that the number of iterations to converge depends weakly on the frequency and the number of subdomains. Fig. \ref{graph:complexity} shows the empirical scaling for one global GMRES iteration. We can observe that both methods have the same asymptotic runtime, but with different constants.

%The proposed numerical method decomposes the computational domain into $L$ layers. The method reveals its true potential when the number of layers $L$ obeys $1 \ll L \ll N$.
%
%Table \ref{table:complexity} gathers the observed complexities, in which $\gamma$ depends on the scaling of the frequency $\omega$ with respect to $n$. For example, if $\omega \sim n^{1/2}$, the regime in which the 5-point stencil discretization is expected to be accurate, then empirically we obtain $\gamma = 5/8$. In this situation, it is possible to scale the number of subdomains with the number of unknowns following $L \sim N^{3/16}$, leading to an overall online complexity of $\cO(N^{13/16})$.
%
%\begin{table}
%\begin{center}
%\begin{tabular}{|c|c|c|c|c|}
%\hline
%Step 		&  Sequential  & Zero-comm parallel    & Number of processors & Communication cost \\
%\hline
%offline 	& - 		& $\cO(\left( N/L \right)^{3/2} )$  		& $ N^{1/2}L$  & $ \cO(\alpha L + \beta NL)$  \\
%\hline
%online			& $\cO(N^{\gamma} L) $  	& $\cO\left( N\log(N) /L \right)$ & $L$   &  $ \cO(\alpha L + \beta N^{1/2}L)$   \\
%\hline
%\end{tabular}
%\caption{Time complexities for the different stages of the solver. The parameter $\gamma$ is \emph{strictly} less than one; its value depends on the scaling of $\omega$ vs. $N$. Here $\alpha$ is the latency and $\beta$ is the inverse bandwidth.
%}\label{table:complexity}
%\vspace{-0.5cm}
%\end{center}
%\end{table}

\section{Discussion}

We presented an extension to \cite{ZepedaDemanet:the_method_of_polarized_traces}, with improved asymptotic runtimes in a distributed memory environment. The method has sublinear runtime even in the presence of rough media of geophysical interests. Moreover, its performance is completely agnostic to the source.

This algorithm is of especial interest in the context of 4D full wave form inversion, continuum reservoir monitoring, and local inversion. If the update to the model is localized, most of precomputations can be reused. The only extra cost is the refactorization and computation of the Green's functions in the cells with a non null intersection with the support of the update, reducing sharply the computational cost.

We point out that this approach can be further parallelized using distributed algebra libraries. Moreover, the sweeps can be pipelined to maintain a constant load among all the nodes.

\section{Acknowledgments}

The authors thank Total SA for support. LD is also supported by AFOSR, ONR, and NSF.

\bibliographystyle{seg}  % style file is seg.bst
\bibliography{GRF_integral_formulations.bib}

\end{document}